%% file: 0FPS-Coherence-FINAL.tex
\documentclass[12pt]{amsart}

\usepackage{amssymb,amscd,amsthm,mathrsfs}
\usepackage[ansinew]{inputenc}
\usepackage[centertags]{amsmath}
\usepackage[all]{xypic}
\usepackage{graphicx}
\usepackage{color}

\newcommand{\R}{\mbox{$\mathbb{R}$}}

\newcommand{\T}{\mbox{$\mathbb{T}$}}

 \def\RR{{\mathbb R}}  \def\TT{{\mathbb T}}
   
 \def\ZZ{{\mathbb Z}}

\def\cC{\mathcal{C}}    \def\cU{\mathcal{U}}
    
\def\cE{\mathcal{E}}    \def\cW{\mathcal{W}}
\def\cF{\mathcal{F}}

\newtheorem*{teo*}{Theorem}
\newtheorem*{mainteo}{Main Theorem}
\newtheorem*{teoA}{Theorem A}
\newtheorem*{teoB}{Theorem B}
\newtheorem*{corC}{Corollary C}
\newtheorem{teo}{Theorem}[section]

\newtheorem{quest}{Question}
\newtheorem{cor}[teo]{Corollary}
\newtheorem{af}[teo]{Claim}
\newtheorem{lema}[teo]{Lemma}
\newtheorem{prop}[teo]{Proposition}

\newcommand{\bi}{\begin{itemize}}
\newcommand{\ei}{\end{itemize}}

\theoremstyle{definition}
\newtheorem{defi}{Definition}
\newtheorem{rem}[teo]{Remark}

\newcommand{\demo}[1]{\vspace{.05in}{\sc\noindent Proof #1.}}
\newcommand{\dem}{\vspace{.05in}{\sc\noindent Proof.\,\,}}
\newcommand{\lqqd}{\par\hfill {$\Box$} \vspace*{.05in}}
\newcommand{\finobs}{\par\hfill{$\diamondsuit$} \vspace*{.05in}}

\newcommand{\eps}{\varepsilon}

\newcommand{\en}{\subset}

\newcommand{\PH}{\mathsf{PH}}

\usepackage{ulem}

\newcommand{\comment}[1]{}

\author[T. Fisher]{Todd Fisher}
\address{Department of Mathematics, Brigham Young University,
Provo, UT 84602}
\urladdr{math.byu.edu/$\sim$tfisher/}
\email{tfisher@math.byu.edu}

\author[R. Potrie]{Rafael Potrie}
\address{CMAT, Facultad de Ciencias, Universidad de la Rep\'ublica, Uruguay}
\urladdr{www.cmat.edu.uy/$\sim$rpotrie}
\email{rpotrie@cmat.edu.uy}

\author[M. Sambarino]{Mart\'in Sambarino}
\address{CMAT, Facultad de Ciencias, Universidad de la Rep\'ublica, Uruguay}
\email{samba@cmat.edu.uy}

\title[Dynamical coherence]{Dynamical coherence of partially hyperbolic diffeomorphisms of tori isotopic to Anosov}
\thanks{T.F. was partially supported by the Simons Foundation grant \# 239708. R.P. and M.S. were partially supported by CSIC group 618 and Balzan's research project of J.Palis. R.P. was also partially supported by FCE-3-2011-1-6749. }

\begin{document}

\begin{abstract}
We show that partially hyperbolic diffeomorphisms of $d$-dimensional tori isotopic to an Anosov diffeomorphism, where the isotopy is contained in the set of partially hyperbolic diffeomorphisms, are dynamically coherent. Moreover, we show a \textit{global stability result}, i.e. every partially hyperbolic diffeomorphism as above is \textit{leaf-conjugate} to the linear one.  As a consequence, we obtain intrinsic ergodicity  and measure equivalence for partially hyperbolic diffeomorphisms with one-dimensional center direction that are isotopic to Anosov diffeomorphisms through such a path.


\bigskip


\medskip

\end{abstract}

\maketitle

\section{Introduction}\label{SectionIntroduccion}

A fundamental problem in dynamical systems is classifying dynamical phenomena and describing the spaces that support these actions.  By the 1970s there was a good classification of smooth systems that are uniformly hyperbolic.

This is seen in the well known Franks-Manning classification result of Anosov diffeomorphisms of tori. This result provides a global classification of Anosov diffeomorphisms on tori (or more generally infranilmanifolds) up to topological conjugacy:  any Anosov diffeomorphism of $\TT^d$ is topologically conjugate to its linear part. The proof uses the structure of the invariant foliation as a key tool in obtaining such a classification. Such a result is sometimes referred to as a \textit{global stability} result since it provides classification beyond small perturbations of the system (which is referred to as \textit{structural stability}).

For Anosov flows there are also some global stability results: let us mention for example a result of Ghys \cite{Ghy} (that overlaps with some related results of Gromov \cite{Gromov}) which states that if $\phi_t$ is an Anosov flow in a 3-manifold which is a circle bundle over a surface then $\phi_t$ is \textit{orbit equivalent} to the geodesic flow in a surface of negative curvature.  In the case of flows, orbit equivalence is the natural extension of topological conjugacy as it is well known that topological conjugacy is too strong an equivalence  for local stability results.

Recently, there is a great deal of interest in understanding the dynamical properties of partially hyperbolic diffeomorphisms, precise definitions are given in Section~\ref{s.precise}. For partially hyperbolic diffeomorphisms the natural equivalence relation is given by \textit{leaf-conjugacy} as introduced in \cite{HPS} where a local stability result is provided (under some technical hypotheses). In dimension 3 some global stability results have been obtained (see \cite{Hammerlindl,HP}). These involve studying integrability of the center bundle since the notion of equivalence up to leaf conjugacy relies on the existence of a center foliation, sometimes referred to as \textit{dynamical coherence}. In dimension 3 the center bundle is one-dimensional, a hypothesis that provides some starting point for studying integrability.

In this paper we consider partially hyperbolic diffeomorphisms of the $d$-torus isotopic to Anosov diffeomorphisms with no restriction on the dimension of the center bundle. The main result provides integrability of the center bundle as well as a global stability statement in the case the partially hyperbolic diffeomorphism can be connected to the linear Anosov diffeomorphism by a path of partially hyperbolic diffeomorphisms. The known techniques of working with codimension one foliations are no longer available and we must trade those by dynamical-geometrical properties of the foliations related with the existence of global semiconjugacies to the linear representative. We start with an informal presentation of our results followed by a more precise formulation.

\subsection{Dynamical coherence}\label{s.coherence}

It is well known that the stable and unstable bundles of an Anosov diffeomorphism are integrable.   This extends to the stable and unstable bundles of a partially hyperbolic diffeomorphism~\cite{HPS}, but the integrability of the center bundle is a subtler issue, see for instance~\cite{BuW}.  When the center bundle is integrable the partially hyperbolic diffeomorphism is {\it dynamically coherent}.

\begin{mainteo}\label{t.main} If $f: \TT^d \to \TT^d$ is a partially hyperbolic diffeomorphism that is isotopic to a linear Anosov automorphism along a path of partially hyperbolic diffeomorphisms, then $f$ is dynamically coherent.
\end{mainteo}

We establish dynamical coherence without the usual restrictions on the dimension of the center bundle, the strength of the domination, or the geometry of the strong foliations in the universal cover.
This is one of the first results on dynamical coherence without restriction on the dimension of the center bundle that holds in ``large'' open sets (whole connected components of partially hyperbolic diffeomorphisms) and the center fibers are noncompact.

The techniques we introduce also allow us to show plaque-expansiveness as defined in \cite{HPS} which in turn give a global stability result in this setting. See section \ref{Section-LeafConj}.

\subsection{Maximizing measures}

Another motivation for this paper grew from an attempt to extend the results of \cite{BFSV}. In that paper it is shown that a well known example of partially hyperbolic diffeomorphism, known as Ma\~n\'e's example (see \cite{ManheContributions} or \cite{BDV} Chapter 7), has a unique measure of maximal entropy. In fact, it is shown there that using the measure of maximal entropy that Ma\~n\'e's example as a measure preserving transformation is isomorphic to the measure preserving transformation given by a linear Anosov automorphism of $\TT^3$ and Haar measure.

This result was extended in~\cite{BF} to certain diffeomorphisms that are $C^0$ close to hyperbolic toral automorphisms, but not partially hyperbolic.  In this case the diffeomorphisms satisfy a weak version of hyperbolicity called a dominated splitting.

A further extension was obtained by Ures~\cite{Ures} to all \textit{absolutely} partially hyperbolic diffeomorphisms of $\TT^3$ isotopic to Anosov as well as other higher dimensional cases under the further assumption of quasi-isometry of the strong foliations (in order to be able to use results of \cite{Brin,Hammerlindl}). For $\TT^3$, under the assumption of pointwise partial hyperbolicity, this result can be weakened to cover all (not necessarily absolute) partially hyperbolic diffeomorphisms of $\TT^3$ isotopic to Anosov thanks to the results in \cite{Pot}, see \cite{HP} section 6.1.
%

Let us briefly comment on the idea of the proof of the existence and uniqueness of maximal entropy measures for partially hyperbolic diffeomorphisms with one dimensional center isotopic to Anosov. For such diffeomorphisms there always exists a continuous semiconjugacy to their linear part, and the main point in the proof consists in showing the following properties:

\begin{itemize}
\item The fibers of the semiconjugacy are connected arcs of bounded length (and thus carry no entropy).
\item  The image of the set of points on which the semiconjugacy is $1$ to $1$ has total Lebesgue measure in $\TT^d$.
\end{itemize}

These two results together with properties of topological and measure theoretic entropy give the desired result (see Section \ref{Section-MEM} for more details). The main point is to obtain dynamical coherence and use the fact that fibers of the semiconjugacy are contained in center manifolds, this is to be expected since one expects the semiconjugacy to be injective along strong manifolds. This is why in \cite{Ures} the hypotheses of quasi-isometry and absolute partial hyperbolicity are used.

We prove that partially hyperbolic diffeomorphisms (not necessarily absolute) which are isotopic to the linear Anosov diffeomorphisms along a path of partially hyperbolic diffeomorphisms with one-dimensional center bundle have a unique measure of maximal entropy, 
see Corollary C. We remark that even for absolute partially hyperbolic diffeomorphisms this result was not known without further hypotheses on the geometry of the strong foliations.

We remark that in~\cite{NY, BF} it was shown that there are systems with a unique measure of maximal entropy and whose topological entropy is $C^1$ locally constant even if the center bundles have dimension 2.  In~\cite{NY} the situation is a partially hyperbolic diffeomorphism that is dynamically coherent with 2-dimensional center fibers, and in~\cite{BF} there are two transverse foliations each 2-dimensional and tangent to the dominated splitting.  In both of these cases the diffeomorphisms can be chosen to be isotopic to Anosov. Moreover, in~\cite{BFSV} it shown that there are partially hyperbolic diffeomorphisms isotopic to Anosov (through a path of partially hyperbolic ones) having bidimensional center and having a unique measure of maximal entropy (and whose topological entropy is also $C^1$ locally constant). This example can be extended to higher dimensional center. Thus,  another reason to establish dynamical coherence in the Main Theorem is that under certain additional hypotheses one may be able
to establish there is a unique measure of maximal entropy and constant topological entropy for systems isotopic to Anosov diffeomorphisms without the restriction of the center bundle being 1-dimensional (although of course one cannot expect that this holds in the entire connected component in this case).

\subsection{Precise Setting}\label{s.precise}

We say that $f: M \to M$ is \textit{partially hyperbolic} if there exists a $Df$-invariant splitting $TM = E^{ss}_f \oplus E^c_f \oplus E^{uu}_f$ such that there exists $N>0$ and $\lambda>1$ verifying that for every  $x\in M$ and unit vectors $v^{\sigma}\in E^\sigma_f(x)$ ($\sigma= ss,c,uu$) we have

\begin{itemize}
\item[(i)] $\lambda \|Df_x^N v^{ss}\| <  \|Df_x^N v^c \| < \lambda^{-1} \| Df_x^N v^{uu} \|, $ and
\item[(ii)] $ \|Df_x^N v^{ss}\| < \lambda^{-1} < \lambda <  \| Df_x^N v^{uu} \|$.
\end{itemize}

We will assume throughout that $N=1$ due to results in \cite{Gourmelon}. We remark that the bundles can be trivial.

The definition we have used of partial hyperbolicity is the weakest one appearing in the literature. It is sometimes referred to as \textit{pointwise} partial hyperbolicity as opposed to \textit{absolute} partial hyperbolicity
\footnote{For absolute partial hyperbolicity it is required that the inequalities hold for unit vectors that may belong to the bundles of different points.}.
The absolute partial hyperbolicity sometimes simplifies proofs of dynamical coherence (see \cite{Brin}) but is quite artificial as it does not capture the real nature of domination (this becomes clear for example when more bundles are involved). We remark that there are different results in the study of dynamical coherence depending on the definition used, see \cite{BBI2,HHU,Pot}.

We denote
$$\PH(\TT^d) = \{ f: \TT^d \to \TT^d  \text{ partially hyperbolic} \}.$$

Let $A \in SL(d,\ZZ)$ be a linear Anosov automorphism admitting a dominated splitting of the form $E^{ss}_A \oplus E^{ws}_A \oplus E^{wu}_A \oplus E^{uu}_A$. We denote as $E^s_A=E^{ss}_A \oplus E^{ws}_A$, $E^c_A= E^{ws}_A \oplus E^{wu}_A$ and  $E^u_A = E^{wu}_A \oplus E^{uu}_A$. There may be many possibilities  for the dimensions of $E^{ss}_A$ and $E^{ws}_A$ (respectively for $E^{wu}_A$  and $E^{uu}_A$).

We consider $\PH_{A,s,u}(\TT^d) \en \PH(\TT^d)$ the subset of those which are isotopic to $A$ and whose splitting verifies $\dim E^{ss}_f = \dim E^{ss}_A=s$ and $\dim E^{uu}_f = \dim E^{uu}_A=u$. In order to simplify notation we will denote $\PH_{A,s,u}(\TT^d)$ as $\PH_A(\TT^d)$ leaving the dimensions of $E^{\sigma}_A$ ($\sigma=ss,uu$) implicit from the context (we will leave them fixed throughout the paper).

For $X\subset\TT^d$ we let $\widetilde X$ denote the lift of $X$ to $\RR^d$.  Similarly, for $f:\TT^d\rightarrow \TT^d$ a diffeomorphism we let $\widetilde{f}:\RR^d\rightarrow \RR^d$ be the lift of $f$.

Given $f \in \PH_A(\TT^d)$ we know from~\cite{Fr} there exists $H_f: \RR^d \to \RR^d$ a continuous and surjective map such that
$$ A \circ H_f = H_f \circ \widetilde f .$$

Moreover, $H_f(x + \gamma) = H_f(x) + \gamma$ for every $\gamma \in \ZZ^d$.

\begin{rem}\label{Remark-ContinuousVariationSemiconj} The map $H$ varies continuously with $f$ in the $C^0$-topology. This is a general fact which does not require $f$ to be partially hyperbolic. This means that given $\eps>0$ there exists a neighborhood $\cU$ of $f$ in the $C^0$-topology such that $d(H_f(x),H_g(x))< \eps$ for every $x\in \RR^d$ and $g\in \cU$.
\finobs
\end{rem}

We say that $f$ is \textit{dynamically coherent} if there exist $f$-invariant foliation $\cW^{cs}_f$ and $\cW^{cu}_f$ tangent respectively to $E^s_f\oplus E^c_f$ and $E^c_f \oplus E^u_f$ (and hence there exists an invariant center foliation $\cW^c_f$ tangent to $E^c_f$).

We say that a dynamically coherent $f \in \PH_A(\TT^d)$ is \textit{center-fibered} if $H_f^{-1}(E^c_A + H_f(x))= \widetilde \cW^c_f(x)$. This means that by the semiconjugacy $H_f$ different leaves of the center foliation map surjectively to different translates of $E^c_A$.

We denote $\PH_A^0(\TT^d)$ to be the connected components of $\PH_A(\TT^d)$ containing a dynamically coherent and center-fibered partially hyperbolic diffeomorphism. Notice that the linear Anosov diffeomorphism $A$ is center-fibered so that $\PH_A^0(\TT^d)$ is a non-empty open set with at least one connected component. Notice also that the space of Anosov diffeomorphisms may not be connected~\cite{FG}, so that the set $\PH_A^0(\TT^d)$ is potentially larger than the connected component containing $A$. Note also that in \cite{FG} the construction is based by showing that there is an Anosov diffeomorphisms which is conjugated to $A$ by a diffeomorphism isotopic to the identity (but not diffeotopic) and therefore this Anosov diffeomorphism is dynamically coherent and center fibered.

At the moment, we do not know the answer to the following questions which we believe to have affirmative answers and would improve our results considerably.

\begin{quest}
Is every partially hyperbolic diffeomorphism in $\PH_A(\TT^d)$ in the connected component of a dynamically coherent and center-fibered one? In other words, does $\PH_A^0(\TT^d)= \PH_A(\TT^d)$?
\end{quest}

In fact, the hypothesis of being center-fibered is crucial to our proofs but it is not clear whether it follows from dynamical coherence or not.

\begin{quest}
Is there an example of a partially hyperbolic diffeomorphism in $\PH_A(\TT^d)$ such that it is dynamically coherent but not center-fibered?
\end{quest}

An affirmative answer to the first question would imply an affirmative answer to the second. However, in view of the results of \cite{FG} it seems clear that (if it admits a positive answer) the first question is much harder in principle.

\subsection{Precise Statement of results}

\begin{teoA}\label{Teo-Main} Every $f \in \PH_A^0(\TT^d)$ is dynamically coherent and center fibered.
\end{teoA}

We prove some intermediary results in more generality. Also, the theorem can be applied even in the case where $E^{ss}_A$ or $E^{uu}_A$ are zero dimensional (if both are trivial, the theorem itself is trivial).

 We also  provide a global stability result in this context by showing \textit{leaf conjugacy} and thus improving previous results on the case of one dimensional center (see \cite{Hammerlindl, HP}). Recall that two dynamically coherent partially hyperbolic diffeomorphisms $f,g:M\to M$ are said to be \textit{leaf conjugate} if there exists a homeomorphism $h:M\to M$ such that $h(\cW_f^c(x))=\cW_g^c(h(x))$ and
 $h\circ f(\cW_f^c(x))=\cW_g^c(g\circ h(x)).$

 \begin{teoB}\label{t.leafconjugacy}
 Any two diffeomorphisms in the same connected component of $PH_A^0(\T^d)$ are leaf conjugate. In particular any diffeomorphism of $PH^0_A$ in the same connected component of $A$ is leaf conjugate to $A.$
 \end{teoB}

We also investigate the existence of measures of maximal entropy and we deduce the following consequence:

\begin{corC}\label{Cor-Main} If $\dim E^c_f=1$ then there exists a unique maximal
entropy measure with equal entropy to the linear part.  
\end{corC}

See Section \ref{Section-MEM} for definitions and the proof of the Corollary.


It is possible that our results can be applied in the case of partially hyperbolic diffeomorphisms isotopic to Anosov diffeomorphisms in nilmanifolds. However, this has to be done with some care since even the initial Anosov diffeomorphism may not be dynamical coherent (see \cite{BuW} for possible problems). It may then be the case that every partially hyperbolic diffeomorphism isotopic to such Anosov through partially hyperbolic diffeomorphisms will not be dynamically coherent (extending a construction announced by Gourmelon \cite{BuW}), but we have not checked this in detail. 

It is also possible that our techniques shed light in studying the case of partially hyperbolic diffeomorphisms of $\TT^d$ isotopic to linear partially hyperbolic automorphisms even if these are not Anosov. This is because there are some types of semiconjugacies when the linear part is partially hyperbolic, and under some (possibly more restrictive) hypotheses one expects that our techniques could be adapted to that case. Notice that the non-dynamical coherent examples given by \cite{HHU} are not isotopic to their linear representative through partially hyperbolic systems.

\textbf{Organization of the paper:} In the next section we provide some basic definitions and preliminary results. In sections \ref{s.proper} and \ref{s.openclose} we state the main property and prove that is an open and closed property in $PH_A^0(\T^d)$ which is fundamental to prove our results. In section \ref{s.dyncoherence} we give sufficient conditions to have dynamically coherence. Theorem 
A is proved in section \ref{s.teomain} and Theorem B 
is proved in Section \ref{Section-LeafConj} together with other results concerning quasi isometric foliations. Section \ref{Section-MEM} is devoted to the study of measures of maximal entropy and to prove Corollary C. 

\section{Definitions and preliminaries:}

\subsection{First remarks}
For $f\in\mathrm{PH}(\TT^d)$
there exist $f$-invariant foliations $\cW^{ss}_f$ and $\cW^{uu}_f$ tangent to $E^{ss}_f$ and $E^{uu}_f$ respectively that we call the \textit{strong foliations}.  We let $\widetilde \cW^\sigma(x)$ denote the associated $\sigma$ foliation for $\widetilde{f}$ (where $\sigma= ss, uu$ or when they exist $\sigma=cs,cu,c$).

In general, we have $H_f(\widetilde \cW^{uu}_f(x)) \en E^{wu}_A \oplus E^{uu}_A +
H_f(x)$. Similarly, for $\widetilde \cW^{ss}_f$ we have $H_f(\widetilde \cW^{ss}_f(x)) \en E^{ws}_A \oplus E^{ss}_A +
H_f(x)$.

We now introduce some notation.  Let
$$B^\sigma_R(x) = B_R(x) \cap (E^\sigma_A + x) $$
for $\sigma = ss, uu, c, s,u, ws,wu$ where $B_R(x)$ is the ball of
raius $R$ centered at $x.$ For $f \in \PH(\TT^d)$ we let
$$ D^{\sigma}_{R,f} (x) = \{ y \in \widetilde \cW^\sigma(x) \ : \ d_{\cW^\sigma}(x, y)
< R \} $$
where  $d_{\cW^\sigma}(\cdot, \cdot)$ denotes the metric inside the leaves induced by restricting the metric of $\RR^d$ to a Riemannian metric in the leaves. Sometimes we will denote $d_{\cW^\sigma}$ as $d_{\sigma}$.

From the continuous variation on compact parts of the strong manifolds one has the following classical result \cite{HPS}.

\begin{prop}\label{Prop-ContinuousVariation} For every $R>0$ and $\eps>0$ there exists $\cU$ a $C^1$-neighborhood of $f$ and $\delta>0$ such that for every $g \in \cU$ and every $x,y\in \RR^d$ with $d(x,y)<\delta$ one has

$$ d_{C^1} (D^\sigma_{R,g}(x), D^\sigma_{R,f}(y)) < \eps $$

\noindent for $\sigma= ss, uu$.
\end{prop}

\begin{rem}\label{Remark-PHConstants} For $f\in \PH(\TT^d)$ there exist constants $1< \lambda_f < \Delta_f$ such that in a $C^1$-neighborhood $\cU$ of $f$ we have
$$ D^{uu}_{(\lambda_f R),g}(\widetilde g(x)) \en \widetilde g (D^{uu}_{R,g}(x)) \en D^{uu}_{(\Delta_f R), x} (\widetilde g(x)) $$
for every $g \in \cU$, $x \in \RR^d$  and $R>0$.
A similar result holds for $D^{ss}$ by applying $\widetilde g^{-1}$. This follows from the fact that the derivative of $f$ restricted to the unstable bundle is always larger than $\lambda_f$ and the global derivative of $f$ is smaller than $\Delta_f$ (from compactness). Therefore, one can also show for $g$ $C^1$-close to $f$ that one has the same estimates for the derivative of $g$ in any vector lying in a small cone around the unstable direction of $f$, so that the estimates hold for disks tangent to a cone close to the unstable direction. \finobs
\end{rem}

\subsection{Strong Almost Dynamical Coherence}

The following definitions are motivated by the ones introduced in \cite{Pot} but slightly adapted to our needs.

\begin{defi}[Almost parallel foliations] Let $\cF_1$ and $\cF_2$ be foliations of $\TT^d$.  Then they are \textit{almost parallel} if
there exists $R>0$ such that for every $x\in \RR^d$ there exists $y_1$ and $y_2$ such that:
\begin{itemize}
\item $\widetilde \cF_1(x) \en B_R (\widetilde \cF_2(y_1))$ and $\widetilde
\cF_2(y_1) \en B_R(\widetilde \cF_1(x))$, and 
\item $\widetilde \cF_2(x) \en
B_R(\widetilde \cF_1(y_2))$ and $\widetilde \cF_1(y_2) \en B_R(\widetilde
\cF_2(x)).$
\end{itemize}\finobs
\end{defi}

Being almost parallel is an equivalence relation (see \cite{HP} Appendix B). Notice that the condition can be stated in terms of Hausdorff distance by saying that for every $x\in \RR^d$ there exists $y_1$ and $y_2$ such that the Hausdorff distance between $\widetilde \cF_1(x)$ and $\widetilde \cF_2(y_1)$ is smaller than $R$ and the Hausdorff distance between $\widetilde \cF_2(x)$ and $\widetilde \cF_1(y_2)$ is smaller than $R$.

\begin{defi}[Strong Almost Dynamical Coherence] Let $f \in \PH_A(\TT^d)$ we say it is \textit{strongly almost dynamically
coherent} (SADC) if there exists foliations $\cF^{cs}$ and $\cF^{cu}$ (not necessarily invariant) which are respectively
transverse to $E^{uu}_f$ and $E^{ss}_f$ and are almost parallel to the foliations $E^{ss}_A \oplus E^c_A$ and $E^c_A \oplus E^{uu}_A$ respectively. \finobs
\end{defi}

The next result is proved in \cite[Proposition 4.5]{Pot}.

\begin{prop}\label{Prop-SADCopenclosed} Being SADC is an open and closed property in $\PH_A(\TT^d)$. In particular, every $f \in \PH_A^0(\TT^d)$ verifies this property.
\end{prop}

The idea of the proof is that open is trivial since the same foliation works by the continuous variation of the $E^{ss}$ and $E^{uu}$ bundles.
If $f_n \to f$ one can choose $n$ large enough so that the bundles are close. By choosing the foliation $\cF^{cs}_n$ for $f_n$ and
iterating backwards by $f_n$ a finite number of times one gets a foliations which works for $f$. Notice that since $f_n$ is isotopic to $A$ it fixes the class of foliations almost parallel to any $A$-invariant hyperplane.

\section{$\sigma$-properness}\label{s.proper}

We define $\Pi^{\sigma}_x$ to be the projection of $\RR^d$ onto $E^\sigma_A + x$ along the complementary subbundles of $A$, we will usually omit the subindex $x$. Let $H^{\sigma}_f:= \Pi^{\sigma}\circ H_f$.

\begin{defi}[$\sigma$-properness] For $\sigma=ss,uu$ we say that $f\in \PH_A(\TT^d)$ is $\sigma$-\textit{proper} if the map $H^\sigma_f|_{\widetilde \cW^\sigma}$ is (uniformly) proper. More precisely, for every $R>0$ there exists $R'>0$ such that, for every $x\in \RR^d$ we have that \footnote{This can also be expressed as: $(H_f^\sigma)^{-1}(B^\sigma_R (H_f(x))) \cap \widetilde \cW^\sigma_f(x) \subset D^\sigma_{R',f}(x)$.}
$y\in \widetilde \cW^{\sigma}(x)$ and $d(H^{\sigma}_f(x), H^{\sigma}_f(y))<R$ implies $d_{\sigma}(x,y)<R'$.
\finobs
\end{defi}

\begin{lema}\label{Lema-AlcanzaUno} Assume $f \in \PH_A(\TT^d)$  such that there exists $R_1>0$ verifying that for every $x\in \RR^d$ we have $y \in \widetilde \cW^{\sigma}(x)$ and $d(H^\sigma_f(x), H^\sigma_f(y)) < 1$ implies $d_\sigma (x,y) < R_1$.
Then $f$ is $\sigma$-proper.
\end{lema}

\dem We consider the case $\sigma=uu$ the other is symmetric. Since $A$ is Anosov and expands uniformly along $E^{uu}_A$ we know that given $R>0$ there exists $N>0$ such that for every $z\in \RR^d$ we have $B^{uu}_R(A^N(z)) \en A^N(B^{uu}_1(z))$ .

Consider $R>0$ and $R' = \Delta_f^N R_1$ with $N$ as defined above and $\Delta_f$ as in Remark \ref{Remark-PHConstants}.
Let
$$y \in (H^\sigma_f)^{-1}(B_R^{uu}(H_f(x))) \cap \widetilde{W}^{\sigma}_f(x).$$
Then, we can see that $\widetilde f^{-N}(y) \in D^{uu}_{R_1}(\widetilde f^{-N}(x))$.
Indeed, since $$H^{uu}_f(y) \in B_R^{uu}(H_f(x))$$ and $A^{-N}(H_f(y)) = H_f(\tilde f^{-N}(y))$
we have that $$\Pi^{uu}(A^{-N}(H_f(y))) \in B^{uu}_1(A^{-N}(x))$$ from how we chose $N$.
Then, from the hypothesis of the Lemma we know that
 $\widetilde f^{-N}(y) \in (H_f^{uu})^{-1}(B^{uu}_1(A^{-N}(x)))$ which is contained in $D^{uu}_{R_1,f}(\widetilde f^{-N}(x))$.

Using Remark \ref{Remark-PHConstants} we deduce that $y \in D^{uu}_{R',f}(x)$ as desired.

\lqqd

In the remainder of this section we will show the equivalence between $\sigma$-properness and the following conditions:

\begin{itemize}
\item[($I^{\sigma}$)] The function $H^\sigma_f$ is injective when restricted to each leaf of $\widetilde \cW^{\sigma}_f$.\\
\item[($S^{\sigma}$)] The function $H^\sigma_f$ is onto $E^\sigma_A+H_f(x)$ when restricted to each leaf of $\widetilde \cW^{\sigma}_f(x)$.
\end{itemize}

\begin{lema}\label{Lema-PropImpliesInjectandSurject} If $f \in \PH_A(\TT^d)$ is $\sigma$-proper, then it verifies both ($I^\sigma$) and ($S^\sigma$).
\end{lema}

\dem First we show the injectivity of $H^\sigma_f$ along leaves of $\widetilde \cW^{\sigma}_f$.
Assume by contradiction that $y$ belongs to the leaf $\widetilde \cW^{\sigma}_f(x)$ of $\widetilde \cW^{\sigma}_f$ and that $H^\sigma_f(x)=H^\sigma_f(y)$ where $y\neq x$. Since $y \neq x$ there exists a $\delta>0$ such that $y \neq D^\sigma_{\delta,f}(x)$. Using Remark \ref{Remark-PHConstants} we know that given $R_1>0$ there exists $N \in \ZZ$ such that $\widetilde f^N(y) \notin D^\sigma_{R_1,f}(\widetilde f^N(x))$.

Consider $R_1$ given by $\sigma$-properness applied to $R=1$. Then, we know that
$$(H^\sigma_f)^{-1}(B^{\sigma}_{1}(H_f(z))) \en D^\sigma_{R_1,f}(z)$$
 for every $z\in \RR^d$.
However, we have
$$(H^\sigma_f)^{-1}(B^{\sigma}_{1} (H_f (\widetilde f^N(x)))$$ contains $\widetilde f^N(y)$, and $\widetilde f^N(y)$ is not contained in $D^\sigma_{R_1,f}(\widetilde f^N(x))$, a contradiction.

Now, we shall show  surjectivity of $H^\sigma_f$ along leaves of $\widetilde \cW^{\sigma}_f$ onto $E^{\sigma}_A$. In the argument we will use the injectivity property established above.

We claim first that injectivity of $H^\sigma_f$ implies that there exists a $\delta>0$ such that
$$H^\sigma_f (\partial D^\sigma_{1,f}(x)) \cap B^{\sigma}_\delta(H_f(x)) = \emptyset. $$
Indeed, otherwise there would exist a pair of sequences $x_n, y_n$ such that $y_n \in \partial D^\sigma_{1,f}(x_n)$ and that
$$H^\sigma_f (y_n) \in B^\sigma_{1/n} (H_f(x_n)).$$
Taking a subsequence and composing with deck transformations we can assume that both sequences converge to points $x,y$. We have that $y \in \partial D^{\sigma}_{1,f}(x)$, in particular $y\neq x$ and we know by continuity of $H_f$ and $\Pi^\sigma$ that $H^\sigma_f(x)=H^\sigma_f(y)$ contradicting injectivity.

From injectivity and Invariance of Domain (see for instance \cite{Hatcher} Theorem 2B.3), we know that for every $z\in \RR^d$ we have $S_z=  H_f^\sigma (\partial D^{\sigma}_{1,f}(x))$ is a $(\dim E^{\sigma}_f -1)$-dimensional sphere embedded in $E^{\sigma}_A + H_f(x)$. Using Jordan's Separation Theorem (\cite{Hatcher} Proposition 2B.1) and the fact that $\dim E^{\sigma}_f = \dim E^{\sigma}_A$ we deduce that $S_z$ separates $E^{\sigma}_A + H_f(x)$ into two components. Moreover, the image by $H^\sigma_f$ of $D^{\sigma}_{1,f}(x)$ is the bounded component and contains $H_f(x)$. From the above remark it also contains $B^{\sigma}_\delta(H_f(x))$.

Now, fix $R>0$, then there exists $N\in \ZZ$ such that
$$B^{\sigma}_R(z) \en A^N(B^\sigma_\delta(A^{-N}(z))).$$
Using the semiconjugacy we see that
$$ B^\sigma_R(H_f(x)) \en  H^\sigma_f (\widetilde f^N(D^\sigma_{1,f}(\widetilde f^{-N}(x)))).$$

Since this holds for any $R$ we know $H^\sigma_f$ verifies ($S^\sigma$) as desired.

\lqqd

\begin{rem}\label{IimpliesS}
Note that in the above proof we have proven that ($I^\sigma$) implies ($S^\sigma$).
\end{rem}

\begin{lema}\label{Lema-ImasSimpliesSigmaProper} If $f \in \PH_A(\TT^d)$ verifies ($I^\sigma$) and ($S^{\sigma}$) then $f$ is $\sigma$-proper.
\end{lema}

\dem The fact that $f$ has properties ($I^\sigma$) and ($S^{\sigma}$) implies that for every $x\in \RR^d$ we know
$$H^\sigma_f : \widetilde \cW^{\sigma}_f(x) \to E^{\sigma}_A + H_f(x)$$
 is a homeomorphism for every $x \in \RR^d$. In particular, we deduce that
 $$(H^\sigma_f)^{-1}(B^{\sigma}_1(H_f(x))) \cap \widetilde \cW^\sigma(x)$$
 is bounded for every $x\in \RR^d$.

Consider the function $\varphi: \RR^d \to \RR$ such that $\varphi(x)$ is the infimum of the values of $R$ such that

$$(H^\sigma_f)^{-1}(B^{\sigma}_1(H_f(x))) \cap \widetilde \cW^\sigma(x) \en D^\sigma_{R,f}(x).$$

\noindent that is to say, the infimum of the values $R$ such that $y\in \widetilde \cW^\sigma(x)$ and $d(H_f(x),H_f(y)) \leq 1$ implies that $d_\sigma(x,y) \leq R$.

From Lemma \ref{Lema-AlcanzaUno} we know that if $\varphi$ is uniformly bounded in $\RR^d$ then $f$ is $\sigma$-proper. Since $\varphi$ is $\ZZ^d$-periodic, it is enough to control its values in a fundamental domain that is compact. Thus, it is enough to show that if $x_n \to x$ then $\limsup \varphi(x_n) \leq \varphi(x)$.

To show this, notice that $H^\sigma_f( D^\sigma_{\varphi(x),f}(x))$ contains $B^\sigma_1(H_f(x))$. Since it is a homeomorphism we deduce that for every $\eps$, there exists $\delta$ such that

$$ B^\sigma_{1+\delta}(H_f(x)) \en H^\sigma_f( D^\sigma_{\varphi(x)+\eps,f}(x)).$$

Using the continuous variation of the $\widetilde \cW^\sigma$ leaves (Proposition \ref{Prop-ContinuousVariation}) and continuity of $H^\sigma_f$ we deduce that for $n$ large enough  that $H^\sigma_f( D^\sigma_{\varphi(x)+\eps,f}(x_n))$ contains $B^\sigma_1(H_f(x_n))$ showing that $\limsup \varphi(x_n) \leq \varphi(x) + \eps$ and this holds for every $\eps>0$.

\lqqd

\section{Dynamical coherence}\label{s.dyncoherence}

We now state a criteria for integrability of the bundles of a partially hyperbolic diffeomorphism.  This criteria  generalizes the one given in \cite{Pot} for dimension 3 (though it requires stronger hypotheses).

We recall that two transverse foliations $\cF_1$ and $\cF_2$ of $\TT^d$ have a \textit{global product structure} if for any two points $x, y \in \RR^d$ the leaves $\widetilde \cF_1(x)$ and $\widetilde \cF_2(y)$ intersect in a unique point.

\begin{teo}\label{Teo-CriteriumCoherence} Assume that $f \in \PH_A(\TT^d)$ verifies the following properties:
\begin{itemize}
\item $f$ is SADC.
\item $f$ is $uu$-proper.
\end{itemize}
Then, the bundle $E^{ss}_f \oplus E^c_f$ is integrable into an $f$-invariant foliation $\cW^{cs}_f$ that verifies
$$ H_f^{-1} ((E^{ss}_A \oplus E^c_A) + H_f(x) )= \widetilde \cW^{cs}_f(x). $$
Moreover, we know $\widetilde \cW^{cs}_f$ has a global product structure with $\widetilde \cW^{uu}_f$.
\end{teo}

\dem We know $\{H^{-1}(E^s_A \oplus E^c_A + y)\}_{y \in \RR^d}$ is an $\widetilde f$-invariant partition of $\RR^d$ that is invariant under deck transformations. This follows as a direct consequence of the semiconjugacy relation and the fact that $H$ is $\ZZ^d$-periodic.
We shall show that under the assumptions of the theorem that $\{H^{-1}(E^s_A \oplus E^c_A + y)\}_{y \in \RR^d}$ is a foliation.

Let $\cF^{cs}$ be a foliation given by the SADC property. Since it is almost parallel to the linear foliation induced by the subspace $E^{ss}_A \oplus E^c_A$ and $H_f$ is a bounded distance from the identity we know $H(\widetilde \cF^{cs}(x))$ is a bounded Hausdorff distance of (a translate of) $E^{ss}_A \oplus E^{c}_A$ for every $x\in \RR^d$.


From the properties ($I^{uu}$) and ($S^{uu}$) we deduce that there is a global product structure between $\widetilde \cF^{cs}$ and $\widetilde \cW^{uu}$. Indeed, consider $x, y \in \RR^d$, we shall first show that $\widetilde \cF^{cs}(x)$ intersects $\widetilde \cW^{uu}(y)$. To do this, consider the set $Q=\RR^d \setminus \widetilde \cF^{cs}(x)$. By a Jordan Separation like result one deduces that the $d-cs-1$-homology of $Q$ is non-trivial where $cs= \dim E^{ss}_A +\dim E^c_A$.  For a proof see for example Lemma 2.1 of \cite{ABP}.

Since $\widetilde \cF^{cs}(x)$ is a bounded Hausdorff distance from $E^{ss}_A \oplus E^c_A$ one deduces that there is a non-trivial cycle of the $d-cs-1$-homology group $H_{d-cs-1}(Q)$ inside $E^{uu}_A$. Choosing this cycle sufficiently far away from $\widetilde \cF^{cs} (x)$, and using properties ($I^{uu}$) and ($S^{uu}$) one deduces the existence of a non-trivial cycle contained in $\widetilde \cW^{uu}_f (y)$. This gives the intersection point (for more details see the proof of Proposition 3.1 in \cite{ABP}).
%
%
%

Now we must prove that the intersection point between $\widetilde \cW^{uu}(x)$ and $\widetilde \cF^{cs}(y)$ is unique. For this, it is enough to show that given a leaf $\widetilde \cF^{cs}(x)$ of $\widetilde \cF^{cs}$ there is no leaf of $\widetilde \cW^{uu}_f$ intersecting $\widetilde \cF^{cs}(x)$ more than once. We will use the following easy facts that follow from the hypotheses we have made on $f$:

\begin{itemize}
\item[(1)] $H_f$ is injective along leafs of $\widetilde \cW^{uu}_f$. Moreover, for every $y\in \RR^d$ we have $H_f(\cW^{uu}_f(y))$ intersects $E^{ss}_A \oplus E^c_A$ in a unique point.
\item[(2)] The image $L = H_f(\widetilde \cF^{cs}(x))$ is contractible and at bounded Hausdorff distance from $E^{ss}_A \oplus E^c_A$.
\end{itemize}

Property (1) and the continuity of $\widetilde \cW^{uu}_f$ and $H_f$ allow us to define a continuous map $\varphi : L \to E^{ss}_A \oplus E^c_A$ that is onto by what we have already proved. Local product structure and property (2) imply that $\varphi$ must be a covering and consequently a homeomorphism. Using again that $H_f$ is injective along $\widetilde \cW^{uu}_f$ we conclude uniqueness of the intersection point as desired.

To finish the proof of the theorem we argue as in Theorem 7.2 of \cite{Pot}. Let us sketch the main points since in this case the proof becomes simpler.
Since $\widetilde \cF^{cs}$ is uniformly transverse to $E^{uu}_f$, there are uniform local product structure boxes in $\RR^d$. Inside each local product structure box, by choosing suitable coordinate systems, one can look at the leaves of the foliations $\widetilde \cF_n = \widetilde f^{-n}(\widetilde \cF^{cs})$ as uniformly bounded graphs from a disk of dimension $cs= \dim E^{ss}_f + \dim E^c_f$ to a disk of dimension $uu = \dim E^{uu}_f$. These family of graphs are precompact in the $C^1$-topology (see for instance \cite{HPS} or \cite[Section 3]{BuW2}). The key point, whose proof is identical as the one of the first claim in the proof of Theorem 7.2 of \cite{Pot} is that the image by $H_f$ of any of these limit graphs (which are $C^1$-manifolds tangent to $E^{ss}_f \oplus E^c_f$) is contained in the corresponding translate of $E^{ss}_A \oplus E^c_A$. Now, using the fact that $H_f$ is injective along strong unstable manifolds, one deduces that such limits are unique, and so the limit graphs form a well
defined foliation with the desired properties (see Theorem 7.2 of \cite{Pot} for more details).

Since the foliation $\widetilde \cW^{cs}_f$ has the same properties of $\widetilde \cF^{cs}$ we get global product structure exactly as above.
\lqqd

A symmetric statement holds for $f$ being $ss$-proper, so we obtain the next corollary.

\begin{cor}\label{Cor-ImplicaCenterFibered} If $f \in \PH_A(\TT^d)$ verifies the SADC property and is both $uu$-proper and $ss$-proper, then $f$ is dynamically coherent and center fibered.
\end{cor}

To prove our main theorem the goal will be to show that having the SADC property and being $\sigma$-proper for $\sigma=uu, ss$ are open and closed properties among partially hyperbolic diffeomorphisms of $\TT^d$ isotopic to linear Anosov automorphisms.

\section{Openness and Closedness of $\sigma$-properness}\label{s.openclose}

In this section we prove that being $\sigma$-proper is an open and closed property among diffeomorphisms in $\PH_A(\TT^d)$ having the SADC property. Without the SADC property it is not hard to show that it is an open property, however, our proof of that it is a closed property uses Theorem \ref{Teo-CriteriumCoherence} so we need  the SADC property (which we already know is open and closed by Proposition \ref{Prop-SADCopenclosed}).

\begin{prop}\label{Prop-OpennesSigmaProper} Being $\sigma$-proper is a $C^1$-open property in $\PH_A(\TT^d)$.
\end{prop}

\dem From Lemma \ref{Lema-AlcanzaUno} it is enough to show that there exists a $C^1$-neighborhood $\cU$ of $f$ such that for each $g\in \cU$ we know there exists $R_1$ such that for every $x\in \RR^d$ we have

$$(H^\sigma_g)^{-1}(B_1^\sigma(H_g(x))) \cap \widetilde \cW^{\sigma}_g(x) \en D^\sigma_{R_1,g}(x).$$

\noindent or equivalently, for every $y \in \widetilde
\cW^\sigma_g(x)$ having $d(H_g^\sigma(x),H_g^\sigma(y))\leq 1$
implies $d_\sigma(x,y) \leq R_1$.

Since $f$ is $\sigma$-proper we know from Lemma \ref{Lema-PropImpliesInjectandSurject} that $H^\sigma_f$ is a homeomorphism from $\widetilde \cW^\sigma_f(x)$ onto $E^\sigma_A + H(x)$. We can choose $R_1$ such that
$$ H^\sigma_f ((D^\sigma_{R_1,f}(x))^c) \cap B^\sigma_2 (H_f(x)) = \emptyset. $$

Let $A^\sigma_{R_1,R_2,g}(x)$ be the annulus $D^\sigma_{R_2,g}(x) \setminus D^\sigma_{R_1,g}(x)$ for any $R_2 > R_1$.  For  $R_2 > \Delta_f R_1$ we have
$$ H^\sigma_f (A^\sigma_{R_1,R_2,f}(x)) \cap B^\sigma_2 (H_f(x)) = \emptyset .$$
Choose $\cU$ a neighborhood of $f$ such that

\begin{itemize}
\item[(i)] the constant $\Delta_f$ holds for every $g\in \cU$ (see Remark \ref{Remark-PHConstants}), and
\item[(ii)] for every $g\in \cU$ we have that $H^\sigma_g (A^\sigma_{R_1,R_2,g}(x)) \cap B^\sigma_1 (H_1(x)) = \emptyset$ (this can be done due to Remark \ref{Remark-ContinuousVariationSemiconj} and Proposition \ref{Prop-ContinuousVariation}).
\end{itemize}

This implies that
$$(H^\sigma_g)^{-1}(B_1^\sigma(H_g(x))) \cap \widetilde \cW^{\sigma}_g(x) \en D^\sigma_{R_1,g}(x).$$
Indeed, otherwise there exists $y \in \widetilde \cW^\sigma_g(x)$ such that $H^\sigma_g (y) \in B_1^\sigma(H_g(x))$ but such that $y \notin D^\sigma_{R_2,g}(x)$.
From the choice of $\Delta_f$ we know that there exists $n\in \ZZ$ such that $\widetilde g^{n}(y) \in A^\sigma_{R_1,R_2}(\widetilde g^n(x))$ (moreover $n>0$ for $\sigma= ss$ and $n< 0$ for $\sigma=uu$) and one knows
$$H^\sigma_g(\widetilde g^n(y)) \in B_1^\sigma(H_g(\widetilde g^n(x)))$$
which contradicts (ii) above.

\lqqd

Notice that the proof shows that the $\sigma$-properness is indeed uniform in the whole neighborhood of $f$.

%

The following is the crucial part in the proof of the main theorem.

\begin{prop}\label{Prop-ClosednessSigmaProper} Being $\sigma$-proper and SADC is a $C^1$-closed property in $\PH_A(\TT^d)$.
\end{prop}

\dem Consider $f_k \to f$ such that $f_k$ are $\sigma$-proper and SADC. From Proposition \ref{Prop-SADCopenclosed} we know that $f$ is also SADC. We will use $k$ instead of $f_k$ in the subscripts to simplify the notation.
\\
Let us assume that $\sigma=uu$. Notice that the diffeomorphisms $f_k$ are in the hypothesis of Theorem \ref{Teo-CriteriumCoherence} so that for every $k>0$ there exist an $f_k$-invariant foliation $\cW^{cs}_k$ tangent to $E^{ss}_k \oplus E^c_k$ which verifies that $\widetilde \cW^{cs}_k(x) = H_k^{-1}((E^{ss}_A \oplus E^c_A) + H_k(x))$, or equivalently:

$$ \text{(}\ast\text{)} \quad H^{uu}_k(x)=H^{uu}_k(y) \ \ \text{ if and only if } \ \ y \in \widetilde \cW^{cs}_k(x).$$

\begin{af}
Given $\eps>0$, there exists $\delta>0$, a cone field $\cC^{uu}$ around $E^{uu}_f$ and $k_0$ such that if $k \geq k_0$ and $D$ is a disk tangent to $\cC^{uu}$ of internal radius larger than $\eps$ and centered at $x$, then $$ B^{uu}_\delta(H_k(x)) \en H^{uu}_k(D). $$
\end{af}

\dem  Consider a finite covering of $\TT^d$ by boxes of local product structure for the bundles of $f$. By choosing them small enough (in particular, smaller than $\eps$) it is possible to assume that the bundles are almost constant in each box (and by changing the metric, also almost orthogonal to each other). By choosing $k_0$ sufficiently large we know that for every $k \geq k_0$ the same boxes are also local product structure boxes for $f_k$.

If $B$ is a box of local product structure we denote by $2B$ and $3B$ the box of double and triple the size, respectively, centered at the same point as $B.$

We can consider the covering small enough and $k_0$ sufficiently large so that for every $k>k_0$ we know

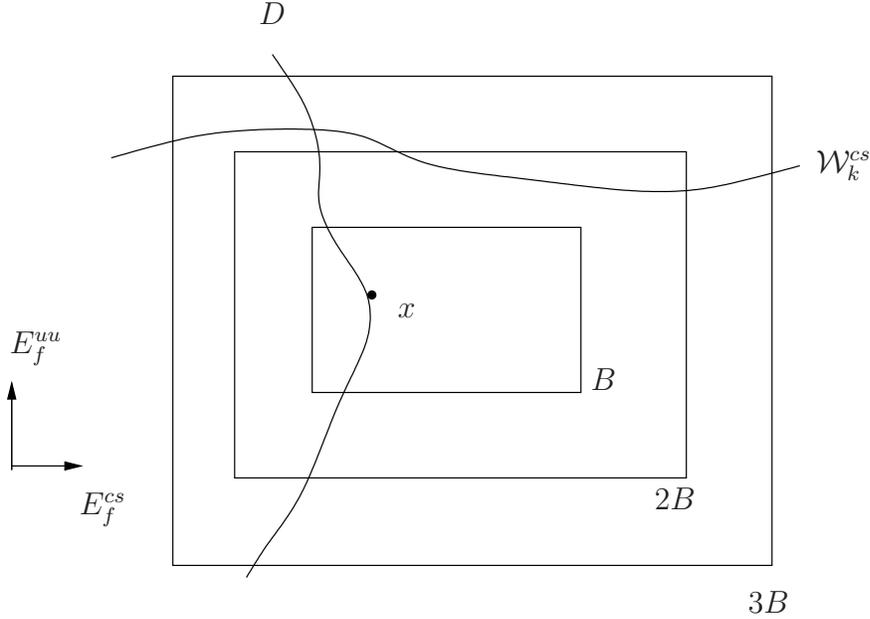
\begin{figure}[ht]\begin{center}
\input{figure.pstex_t}
\caption{\small{The local product structure boxes.}} \label{figure}
\end{center}\end{figure}

\begin{itemize}
\item the boxes $2B$ and $3B$ are also local product structure boxes for all the $f_k$'s in particular
\item for every $B$ of the covering and every disk $D$ tangent to $\cC^{uu}$ of internal radius $\eps$ and centered at a point $x \in B$ we have that $D$ intersects in a (unique) point in $3B$ every center-stable plaque of $\cW^{cs}_k$ which intersects $2B$ (see Figure~\ref{figure}).
\item the previous condition together with ($\ast$) implies that for every disk $D$ tangent to $\cC^{uu}$ of internal radius $\eps$ and centered in a point $x\in B$ one has that $H^{uu}_k(2B) \en H^{uu}_k(D)$.
\end{itemize}

By ($\ast$) above we know that $H_k$ is injective along leaves of $\widetilde \cW^{uu}_k$ so that we have that given a connected component $2B$ of the lift of a local product structure box we have
$$\mathrm{int}( H^{uu}_k (2B))\neq \emptyset. $$

\noindent and that any point $x$ in $B$ is in the interior of $H^{uu}_k(2B)$.

Moreover, since there are finitely many such boxes, we know that there exists a uniform $\delta$ such that $H^{uu}_k(B)$ is at distance $\delta$ from the boundary of $H^{uu}_k(2B)$ independently of the box $B$.

We deduce that every disk $D$ of internal radius $\eps$ centered at a point $x$ and tangent to a small cone around $E^{uu}_f$ verifies that $H^{uu}_k (D)$ contains $B^{uu}_\delta(H_k(x))$ as desired.
\finobs

\begin{af} For any $k$ large enough and $x,y \in \RR^d$ we have that $\widetilde \cW^{uu}_f(x)$ intersects $\widetilde \cW^{cs}_k(y)$.
\end{af}

\dem The previous claim implies that if $k$ is large enough, for every $x, y \in \RR^d$ and we denote $d$ as the distance between  $H^{uu}_k(x)$ and  $ H^{uu}_k(y)$ and let $N_0>\frac{d}{\delta}$, then
\begin{equation}\label{eq.intersection}
 D^{uu}_{N_0\eps, f} (x) \cap \widetilde \cW^{cs}_k(y) \neq \emptyset.
\end{equation}

Indeed, consider the straight segment joining $H^{uu}_k(x)$ with $H^{uu}_k(y)$ in $E^{uu}_A+H_k(x).$ We can cover this segment by $N_0$ balls $B_1,...,B_{N_0}$ of radius $\delta/2$ and such that $B_i\cap B_{i+1}\neq\emptyset.$ Now, $H^{uu}_k(D^{uu}_{\eps, f}(x))$ contains $B_1.$ Thus, $H^{uu}_k(D^{uu}_{2\eps, f}(x))$ contains $B_1\cup B_2$ and inductively we have $H^{uu}_k(D^{uu}_{N_0\eps, f}(x))$ contains $B_1\cup\ldots\cup B_{N_0}$ and $H^{uu}_k(y).$ Using property ($\ast$) above, this implies~\eqref{eq.intersection}.

Therefore, for every $x,y \in \RR^d$ we have that $\widetilde \cW^{uu}_f(x)$ intersects $\widetilde \cW^{cs}_k(y)$.
\finobs

\begin{af} For $k$ large enough the foliations $\widetilde \cW^{uu}_f$ and $\widetilde \cW^{cs}_k$ have a global product structure. Equivalently, the map $H^{uu}_k|_{\widetilde \cW^{uu}_f(x)} : \widetilde \cW^{uu}_f(x) \to E^{uu}_A + H_f(x)$ is a homeomorphism.
\end{af}

\dem By the previous claim, it is enough to show that the
intersection point between $\widetilde \cW^{uu}_f(x)$ and
$\widetilde \cW^{cs}_k(y)$ is unique for any $x,y.$

Since $\widetilde \cW^{uu}_f(x)$ intersects transversally $\widetilde \cW^{cs}_k(y)$ for any $x,y$ and $$H_k(\widetilde \cW^{cs}_k(y))= (E^{ss}_A \oplus E^c_A) + H_k(y)$$ we conclude that $H_k(\widetilde \cW^{uu}_f(x))$ is topologically transversal to $(E^{ss}_A\oplus E^c_A) + H_k(y)$ for any $x,y.$ This implies that $$\Pi^{uu}: H_k(\widetilde\cW^{uu}_f(x))\to E^{uu}_A$$
is a covering map and since $H_k(\widetilde\cW^{uu}_f(x))$ is contractible we know it is one-to-one. Thus, we  have that $H^{uu}_k$ restricted to $\widetilde\cW^{uu}_f(x)$ is a homeomorphism onto $E^{uu}_A$ which implies the desired global product structure.

\finobs

We now return to the proof of the proposition.

We must show (see Lemma \ref{Lema-AlcanzaUno}) that there exists some $R>0$ such that for every $x\in \RR^d$ we have that if $y \in \widetilde \cW^{uu}_f(x)$ and $d(H^{uu}_f(x),H^{uu}_f(y)) \leq 1$ then $d_{uu}(x,y) \leq R$. Equivalently, that

$$ (H^{uu}_f)^{-1}(B^{uu}_1(H_f(x))) \cap \widetilde \cW^{uu}_f(x) \en D^{uu}_{R,f}(x).$$

We will show that for every $x\in \RR^d$ there exists some finite $\psi(x)$ such that

$$(H^{uu}_f)^{-1}(B^{uu}_1(H_f(x))) \cap \widetilde \cW^{uu}_f(x) \en D^{uu}_{\psi(x),f}(x).$$

Then, one can conclude by arguing as in the proof of Lemma \ref{Lema-ImasSimpliesSigmaProper} by considering the infimum $\varphi(x)$ of all possible values of $\psi(x)$ satisfying the property which will be a semicontinuous and periodic function that by a compactness argument is enough to complete the proof.

We know that $d_{C^0}(H_k,H_f)< K_0.$ The previous claim and the fact that $f_k$ is center-fibered implies that $H^{uu}_k$ restricted to $\widetilde\cW^{uu}_f(x)$ is a homeomorphism onto $E^{uu}_A$, so we know that there exists some $R_1>0$ such that
$$H^{uu}_k((D^{uu}_{R_1,f}(x))^c)\cap B^{uu}_{2+2K_0}(H_k(x))=\emptyset$$
and so
$$H^{uu}_f((D^{uu}_{R_1,f}(x))^c)\cap B^{uu}_{1}(H_f(x))=\emptyset.$$
This implies that
$$(H^{uu}_f)^{-1}(B^{uu}_1(H_f(x))) \cap \widetilde \cW^{uu}_f(x) \en D^{uu}_{R_1,f}(x).$$
Setting $\psi(x)=R_1$ we conclude the proof.

\lqqd

\begin{rem} In principle, being $\sigma$-proper could be a closed property too, but in the proof we had to assume also that the sequence (and limit) also had the SADC property. It could be that other proof without the use of that property is possible.
\end{rem}

\section{Proof of Theorem A}\label{s.teomain} 

From our previous results we obtain the following:

\begin{teo}\label{Teo-Main2} If $f \in \PH_A(\TT^d)$ is in the same connected component of a partially hyperbolic $g$ that is $\sigma$-proper (for $\sigma= ss,uu$) and has the SADC property, then $f$ is dynamically coherent and center fibered.
\end{teo}

\dem  Propositions \ref{Prop-OpennesSigmaProper} and \ref{Prop-ClosednessSigmaProper} together with Proposition \ref{Prop-SADCopenclosed} imply that being $\sigma$-proper ($\sigma=ss,uu$) and having the SADC property is an open and closed property in $\PH_A(\TT^d)$.
This implies that every $f$ in the the same connected component of a partially hyperbolic $g$ that is $\sigma$-proper (for $\sigma= ss,uu$) and has the SADC property as in the hypothesis of Corollary \ref{Cor-ImplicaCenterFibered}.
\lqqd

\demo{of Theorem A}
\,
It is enough to show that if $f$ is a partially hyperbolic diffeomorphism in $PH_A^0(\T^d)$ that is  dynamically coherent and center-fibered, then it must be $\sigma$-proper for $\sigma=ss,uu$ and have the SADC property.

This follows from the following remarks:

\begin{itemize}

\item The central stable foliation $\cW^{cs}_f$ is transversal to $E^{uu}_f$ and the central unstable foliation $\cW^{cu}_f$ is transversal to $E^{ss}_f.$

\item Since it is center fibered we know the semiconjugacy is injective along strong stable and unstable manifolds and also that $$H(\cW^{cs}_f(x))\subset (E^{ss}_A\oplus E^c_A) +H(x)$$and
    $$H(\cW^{cu}_f(x))\subset (E^{uu}_A\oplus E^c_A) +H(x).$$

\item Again, since it is center fibered, we know $\Pi^\sigma\circ H$ is injective along strong stable and unstable manifolds. This also implies surjectivity (see Remark \ref{IimpliesS}) and hence we have $\sigma$-properness for $\sigma=ss,uu.$

\item The surjectivity above and the center fibered property implies that
$$H(\cW^{cs}_f(x))= (E^{ss}_A\oplus E^c_A) +H(x)$$and
    $$H(\cW^{cu}_f(x))= (E^{uu}_A\oplus E^c_A) +H(x)$$and from this we easily have the SADC property since $H$ is at bounded distance from the identity.

\end{itemize}
\lqqd


\section{Leaf conjugacy and global stability}\label{Section-LeafConj}

We will now deduce some more additional properties of the systems. We recall that a foliation $\widetilde \cF$ of $\RR^d$ is called \textit{quasi-isometric} if there exist constants $C,D>0$ such that for any pair of points $x,y$ in the same leaf of $\widetilde \cF$ one has

$$ d_{\widetilde \cF}(x,y) \leq C d(x,y) + D $$

\noindent where as before $d_{\widetilde \cF}(\cdot,\cdot)$ denotes the leafwise distance between points and $d(\cdot,\cdot)$ the usual distance in $\RR^d$. We remark that if the foliation $\widetilde \cF$ has $C^1$-leaves, it is possible to change the constants to have $D=0$.

\begin{prop}\label{Proposition-QuasiIsometry} If $f\in \mathrm{PH}^0_A$ is $\sigma$-proper $(\sigma=ss,uu)$ then the foliation $\widetilde \cW^{\sigma}$ is quasi-isometric.
\end{prop}

\dem First we choose a metric on $\R^d$ by declaring $E^{ss}_A, E^c_A$ and $E^{uu}_A$ mutually orthogonal, this metric is equivalent to the usual metric on $\R^d.$ The proof consists of $3$ steps:

\begin{itemize}
\item[(i)] For every $K>0$ there exists $C_K$ such that if $d(x,y)<K$ and $y\in \widetilde \cW^\sigma(x)$ then $d_\sigma(x,y) < C_K d(x,y) $.
\item[(ii)] For every $C_1>0$ there exists $K$ such that for every $x\in \RR^d$ we have that $\widetilde \cW^{\sigma}(x)$ is contained in $B_{K/2}(x) \cup (\cE^\sigma_{C_1}+x)$ where $\cE^\sigma_{C_1}$ is the cone around $E^{\sigma}_A$ of vectors $v = v^\sigma + v^\perp$ satisfying $\|v^\perp \| < C_1 \|v^\sigma \|$ with $v^\sigma \in E^{\sigma}_A$ and $v^\perp \in (E^\sigma_A)^\perp$. Notice that $(E^\sigma_A)^\perp=E^{cs}_A$ if $\sigma=uu$ and $(E^\sigma_A)^\perp=E^{cu}_A$ if $\sigma=ss.$
\item[(iii)] If $y \in \widetilde \cW^\sigma(x)$ one can choose points $x=x_1, \ldots, x_n=y$ in $\widetilde \cW^\sigma(x)$ and $K$ such that $d(x_i,x_{i+1})< K$ and such that $$\sum_{i=1}^{n-1} d(x_i,x_{i+1}) \leq 3 d(x,y).$$
\end{itemize}

Once we have this, putting together properties (i) and (iii) we deduce that

$$ d_\sigma(x,y) \leq \sum d_\sigma (x_i, x_{i+1}) \leq C_K \sum d(x_i,x_{i+1}) < 3C_K d (x,y) $$

\noindent showing quasi-isometry.

We first notice that (i) is a direct consequence of $\sigma$-properness. In fact, if (i) did not hold we would obtain a sequence $x_n,y_n$ of points at distance smaller or equal to $K$ such that $d_\sigma(x_n,y_n) \geq n$. Using $\sigma$-properness we would obtain that $d(H_f(x_n),H_f(y_n)) \to \infty$. On the other hand, since $H_f$ is at bounded distance from the identity and $x_n$ and $y_n$ are at distance smaller than $K$ one gets that $d(H_f(x_n),H_f(y_n))$ must be bounded, a contradiction.

Let us prove (ii). Since $H_f$ is a bounded distance from the identity, to prove (ii) it is enough to show the same property for $H_f(\widetilde \cW^{\sigma}_f(x))$.

Assume that it is not true. Then, there exists a cone $\cE^\sigma$ and we may find sequences (using $\sigma$-properness) $x_n, y_n \in \RR^d$ such that $y_n \in H_f(\widetilde \cW^{\sigma}_f(x_n))$ with $d(y_n,H_f(x_n))\to\infty$ and $y_n\notin \cE^\sigma+H_f(x_n).$ We assume for simplicity that $\sigma=uu,$ the other case is quite similar.

Let $\lambda_c^{-1}=\|A_{/E^{cs}_A}\|$ and let $\lambda_u=\|A^{-1}_{/E^{uu}_A}\|.$ Notice that $\lambda_u/\lambda_c<1.$ Notice first that if $\lambda_c >1$ we know $A_{/E^{cs}_A}$ is contacting so that $H_f : \widetilde \cW^{uu}_f(x) \to E^{u}_A + H_f(x)$ is a homeomorphism and property (ii) is immediate. Also, since $A$ is Anosov we can assume that (maybe by considering an iterate) that $\lambda_c \neq 1$. So, in what follows we shall assume that $\lambda_c < 1$.

Let $\epsilon >0$ and let $m_n=\inf\{m\ge 0:\lambda_c^md(y_n,H_f(x_n))\le \epsilon\}.$ Since
$d(y_n,H_f(x_n))\to\infty$ and $y_n\notin \cE^\sigma+H_f(x_n)$ we have that $m_n\to\infty.$ And we know  $$d(A^{-m_n}(y_n), A^{-m_n}(H_f(x_n)))\ge\lambda_c\epsilon.$$ On the other hand,
$$
\begin{array}{llll}
& d(\Pi^{uu}(A^{-m_n}(y_n)), A^{-m_n}(H_f(x_n)))\\
= & d(A^{-m_n}(\Pi^{uu}(y_n)), A^{-m_n}(H_f(x_n)))\\
\le & \lambda_u^{m_n}\frac{d(y_n,H_f(x_n))}{C_1}\le \left(\frac{\lambda_u}{\lambda_c}\right)^{m_m}\frac{\epsilon}{C_1}\to_{n\to\infty}0.
\end{array}
$$
Now, composing with deck transformation we may assume that $$f^{-m_n}(x_n)\to x$$ and $$A^{-m_n}(y_n)\to y\in H_f(\widetilde \cW^{\sigma}_f(x)),\; y\neq H_f(x).$$ But $\Pi^{uu}(y)=x,$ a contradiction with property ($I^{uu}$) (which follows from $uu$-properness by Lemma \ref{Lema-PropImpliesInjectandSurject}).
Thus,  we obtain that (ii) is verified.

Finally, to prove (iii) we use (ii): We choose $C_1\le 1/2$ and $K$ from (ii) and we define the sequence $x_i$ inductively. First, we impose $x_1=x$. Then, if $d(x_i,y)<K$ we choose $x_{i+1}=y$. Otherwise we pick $x_{i+1}$ as follows. Notice that $d(\Pi^\sigma(y),x_i)\ge \frac{2}{3}K$ and let $z_{i+1}$ be the point in the segment joining $x_i$ and $\Pi^\sigma(y)$ (which is contained in $E^\sigma_A+x_i$) at distance $\frac{2}{3}K$ from $x_i.$

Now, $(\Pi^{\sigma})^{-1}(z_{i+1})\cap (\cE^\sigma +x_i)$ is a disc $D_i$ of radius $\frac{2}{3}C_1 K$ in $(E^\sigma_A)^\perp+z_{i+1}.$ Since $H^\sigma_f$ is homomorphism onto $E^\sigma_A+H_f(x_i)$ when restricted to $\widetilde\cW^\sigma(x_i)$ and $H_f$ is at bounded distance from the identity,  we conclude that $\Pi^\sigma$ is onto $E^\sigma_A+x_i$ when restricted to $\widetilde\cW^\sigma_f(x_i)$.  By (ii)  there is at least one point in $D_i\cap \cW^\sigma_f(x_i).$ We set $x_{i+1}$ to be one of these points.

 We must now show that the process finishes in finitely many steps. Notice that since $y\in \cE^\sigma+x_i$ the straight line segment joining $x_i$ and $y$ intersects $D_i$ and $d(y,D_i)\le d(x_i,y)-\frac{2}{3}K.$ Thus
$$d(x_{i+1},y)\le d(x_i,y)-\frac{2}{3}K +\frac{2}{3}C_1K\le d(x_i,y)-\frac{1}{3}K.$$
So that the process ends in finitely many steps.

Notice also that $d(x_i,x_{i+1})\le K$ and so the above inequality also shows that
$$d(x_{i+1},y)\le d(x_i,y)-\frac{1}{3}d(x_i,x_{i+1}).$$
Therefore, if we have chosen the sequence $x=x_1,x_2,...,x_n=y$ we have by induction that
$$d(x_{n-1},y)\le d(x,y)-\frac{1}{3}\sum_{i=1}^{n-2}d(x_i,x_{i+1})$$ and so
$$\sum_{i=1}^{n-1}d(x_i,x_{i+1})\le 3 d(x,y).$$
\lqqd

When the central dimension is one it is possible  to use the results of \cite{Hammerlindl} to obtain a property called \textit{leaf conjugacy}. This notion is related with the existence of the semiconjugacy but slightly different, it says that there exists a homeomorphism $h: \TT^d \to \TT^d$ which sends center leaves of $f$ to center leaves of the linear Anosov diffeomorphism and conjugates the dynamics modulo the center behavior (see \cite{Hammerlindl} for more details).

The results in~\cite{Hammerlindl} are proved in the absolute partially hyperbolic setting, but in \cite{HP} it is explained which hypothesis should be added in the pointwise case in order to recover his results.

\begin{prop}\label{Cor-ImplicaLeafConjCenter1} Let $f \in \PH_A(\TT^d)$ with $\dim E^c_f=1$ and verifying SADC property and $\sigma$-properness for $\sigma=ss,uu$ then $f$ is leaf conjugate to $A$.
\end{prop}

\dem Theorem 3.2 in \cite{HP} states that the following properties of a dynamically coherent partially hyperbolic diffeomorphism with one dimensional center and isotopic to $A$ guarantee leaf conjugacy.

\begin{itemize}
\item[(i)] The foliations $\widetilde \cW^{\sigma}_f$ ($\sigma=cs,cu$) are almost parallel to the corresponding linear foliations of $A$.

\item[(ii)] The foliations $\widetilde \cW^\sigma_f$ are \textit{asymptotic} to $E^{\sigma}_A$ (i.e. We have that
$$\frac{d(\Pi^\sigma(x),\Pi^\sigma(y))}{d(x,y)} \to 1$$ as $d(x,y) \to \infty$ with $x,y$ in the same leaf of $\widetilde \cW^\sigma$).
\item[(iii)]The foliations $\widetilde \cW^{\sigma}_f$ ($\sigma=ss,uu$) are quasi-isometric.
\end{itemize}

SADC property implies property (i).

It is quite easy to see that using the semiconjugacy with $A$ that conditions ($I^{\sigma}$) and ($S^\sigma$) imply property (ii). Recall that $\sigma$-properness implied properties ($I^\sigma$) and ($S^\sigma$) (Lemma \ref{Lema-PropImpliesInjectandSurject}).

The proof then concludes by applying Proposition \ref{Proposition-QuasiIsometry} to conclude that (iii) is also satisfied.
\lqqd

Using the concept of plaque-expansiveness we are able to prove the previous result without assuming one dimensionality of the center bundle. We remark that in \cite{HamPlaque} it is proved that \textit{absolute partial hyperbolicity} and quasi-isometry implies plaque-expansiveness.

We recall the definition of plaque expansiveness from \cite{HPS}: Let $f: M \to M$ be a dynamically coherent partially hyperbolic diffeomorphism with center foliation $\cW^c$, we say that $f$ is \textit{plaque-expansive} if for every $\eps>0$ sufficiently small the following holds:

\begin{itemize}
\item[--] Let $\{x_n \}_{n\in \ZZ}$ and $\{y_n\}_{n\in \ZZ}$ be two sequences in $M$ such that $d(x_n,y_n) < \eps$ and such that the points $f(x_n)$, $x_{n+1}$ (resp. $f(y_n)$, $y_{n+1}$) belong to the same leaf of $\cW^c$ and $d_{\cW^{c}}(f(x_n), x_{n+1})< \eps$  (resp. $d_{\cW^{c}}(f(y_n), y_{n+1})< \eps$) for every $n$. Then, $x_0$ and $y_0$ belong to the same center leaf and $d_{\cW^{c}}(x_0,y_0)< K \eps$ for $K$ independent of $\eps$.
\end{itemize}

As before, we use $\sigma$-properness to obtain this property:

\begin{prop}\label{prop-PlaqueExpansive}
  Let $f: \TT^d \to \TT^d$ be a partially hyperbolic diffeomorphism isotopic to Anosov such that it is dynamically coherent and center fibered  then the center foliation is plaque expansive.
\end{prop}

\dem Consider $\eps$ small enough so that two points in an $\eps$-neighborhood belong to the same local product structure box (as in the proof of the claim inside the proof of Proposition \ref{Prop-ClosednessSigmaProper}). Note that if we lift this product box to the universal cover (and take a connected component) then we know that a center leaf intersects this box at most in one connected component (otherwise this violates the center fibered property and the injectivity of the semiconjugacy along the strong stable and unstable foliations).

Let $\{x_n\}$ and $\{y_n\}$ be two pseudo-orbits verifying the properties above, that is:

\begin{itemize}
\item[(i)] $d(x_n, y_n) < \eps$ for every $n\in \ZZ$.
\item[(ii)] $d_{\cW^{c}}(f(x_n),x_{n+1})< \eps$ for every $n\in \ZZ$. In particular, $f(x_n)$ and $x_{n+1}$ belong to the same center leaf.
\item[(iii)]$d_{\cW^{c}}(f(y_n),y_{n+1})< \eps$ for every $n\in \ZZ$. In particular, $f(y_n)$ and $y_{n+1}$ belong to the same center leaf.
\end{itemize}

We must prove that this implies that $x_0$ and $y_0$ belong to the same leaf of $\cW^c$ and $d_{\cW^{c}}(x_0,y_0)<K \eps$ for some uniform $K$ which does not depend on $\eps$.

To do this, we will lift the sequences to the universal cover. Choose $\tilde x_0$ and $\tilde y_0$ in $\RR^d$ such that they project respectively to $x_0$ and $y_0$ and such that $d(\tilde x_0,\tilde y_0) < \eps$. Since $\eps$ is small, this already determines uniquely a pair of sequences $\{ \tilde x_n \}$ and $\{ \tilde y_n \}$ satisfying properties analogous to (i), (ii) and (iii).

We now consider the sequences $\{ X_n= H_f(\tilde x_n) \}$  and $\{ Y_n = H_f(\tilde y_n) \}$. Using the expansivity of $A$ one deduces that $X_0$ and $Y_0$ lie in the same leaf of the foliation by translates of $E^{c}_A$. Since $f$ is center-fibered we deduce that both $\tilde x_0$ and $\tilde y_0$ lie in the same leaf of $\tilde \cW^{c}_f$. Moreover, from how we chose the value of $\eps$ and the fact that $d(\tilde x_0, \tilde y_0) < \eps$ one gets the desired property.

\lqqd

As a consequence of this and results of Chapter 7 in \cite{HPS} we obtain our global stability result:

\demo{of Theorem B}
\,
Let $f$ and $g$ be two diffeomorphisms in the same connected component of $PH_A^0(\T^d).$ Consider a path $f_t$ such that $f_0=g$ and $f_1=f$ and such that $f_t$ is in $\PH_A^0(\TT^d)$ for every $t\in [0,1]$.

By the main theorem, $f_t$ is dynamically coherent and center fibered for every $t\in [0,1]$. Thus Proposition \ref{prop-PlaqueExpansive} applies and we know that the center foliation of $f_t$ is plaque-expansive at every $t$. Using the main result of Chapter 7 of \cite{HPS} we deduce that for every $t_0$ there exists $\eps$ such that the diffeomorphism $f_{t_0}$ is leaf-conjugate to every $f_t$ with $t \in (t_0 + \eps, t_0 -\eps)$. By compactness and transitivity of the relation of being leaf conjugate, one deduces  Theorem B. 
\lqqd

\section{Measures of Maximal Entropy}\label{Section-MEM}

The variational principle states that if $f:X\rightarrow X$ is continuous and $X$ is a compact metric space, then $h_{top}(f) = \sup_{\mu} h_\mu(f)$ where $\mu$ varies among all $f$-invariant Borel probability measures, see for instance \cite{Manhe}. It is thus an interesting question to know whether a given system has measures with entropy equal to the topological entropy of the system, and when such measures exist (which are called \textit{measures of maximal entropy}) to know how many of them are there. When there is a unique measure of maximal entropy the system is \textit{intrinsically ergodic}.

Corollary C 
states that every diffeomorphism in $\PH_A^0(\TT^d)$ with one dimensional center is intrinsically ergodic. In \cite{Ures} a similar result is proved under the added assumption of   absolute partially hyperbolic case and under stronger assumptions on the foliations of $f$.

\demo{of Corollary C}
\,
From Theorem A 
and Proposition \ref{Cor-ImplicaLeafConjCenter1}, let $h$ be the semiconjugacy from $f$ to $A$, then for each $x\in\TT^d$ we know that $[x]=h^{-1}h(x)$ is a point or bounded closed interval in the center fiber containing $x$.

The Leddrappier-Walters type arguments in~\cite{BFSV} allow us to conclude that $h_{\mathrm{top}}(f)=h_{\mathrm{top}}(A)$ and that a lift of the Haar measure, $\mu$, for $A$ is a measure of maximal entropy for $f$.

From Lemma 4.1 in~\cite{Ures} we know that
$$\mu\{ x\in\TT^d\, :\, [x]=\{x\}\}=1.$$
Theorem 1.5 in~\cite{BFSV} now applies and we know that $f$ is intrinsically ergodic.  Furthermore, the unique measure of maximal entropy can be seen as the unique lift of Haar measure for $A$ and also as the limit of the measures given by the periodic classes as defined in~\cite{BFSV}.

%

\lqqd

\end{document}

%% file: figure.pstex_t
\begin{picture}(0,0)%
\includegraphics{figure.pstex}%
\end{picture}%
\setlength{\unitlength}{3947sp}%
\begingroup\makeatletter\ifx\SetFigFont\undefined%
\gdef\SetFigFont#1#2#3#4#5{%
  \reset@font\fontsize{#1}{#2pt}%
  \fontfamily{#3}\fontseries{#4}\fontshape{#5}%
  \selectfont}%
\fi\endgroup%
\begin{picture}(5087,3920)(202,-3546)
\put(2664,-1624){\makebox(0,0)[lb]{\smash{{\SetFigFont{12}{14.4}{\rmdefault}{\mddefault}{\updefault}$x$}}}}
\put(3876,-2086){\makebox(0,0)[lb]{\smash{{\SetFigFont{12}{14.4}{\rmdefault}{\mddefault}{\updefault}$B$}}}}
\put(4276,-2836){\makebox(0,0)[lb]{\smash{{\SetFigFont{12}{14.4}{\rmdefault}{\mddefault}{\updefault}$2B$}}}}
\put(4864,-3486){\makebox(0,0)[lb]{\smash{{\SetFigFont{12}{14.4}{\rmdefault}{\mddefault}{\updefault}$3B$}}}}
\put(1789,214){\makebox(0,0)[lb]{\smash{{\SetFigFont{12}{14.4}{\rmdefault}{\mddefault}{\updefault}$D$}}}}
\put(5289,-711){\makebox(0,0)[lb]{\smash{{\SetFigFont{12}{14.4}{\rmdefault}{\mddefault}{\updefault}$\cW^{cs}_k$}}}}
\put(226,-1849){\makebox(0,0)[lb]{\smash{{\SetFigFont{12}{14.4}{\rmdefault}{\mddefault}{\updefault}$E^{uu}_f$}}}}
\put(664,-2861){\makebox(0,0)[lb]{\smash{{\SetFigFont{12}{14.4}{\rmdefault}{\mddefault}{\updefault}$E^{cs}_f$}}}}
\end{picture}%